\documentclass[11pt]{amsart}
\usepackage{pstricks,pst-plot}
\usepackage{mathrsfs}

\parskip=\smallskipamount

\newtheorem{theorem}{Theorem}[section]

\newtheorem{lemma}[theorem]{Lemma}
\newtheorem{corollary}[theorem]{Corollary}

\theoremstyle{definition}

\numberwithin{equation}{section}

\input cyracc.def

\font\tenscr=rsfs10 
\font\sevenscr=rsfs7 
\font\fivescr=rsfs5 
\skewchar\tenscr='177 \skewchar\sevenscr='177 \skewchar\fivescr='177
\newfam\scrfam \textfont\scrfam=\tenscr \scriptfont\scrfam=\sevenscr
\scriptscriptfont\scrfam=\fivescr
\def\scr{\fam\scrfam}

\def\sC{{\scr C}}

\def\sP{{\scr P}}

\def\bB{\mathbb B}\def\bC{\mathbb C}

\def\bM{\mathbb M}\def\bP{\mathbb P}
\def\bR{\mathbb R}\def\bS{\mathbb S}\def\bT{\mathbb T}
\def\bU{\mathbb U}

\def\vp{\varphi}
\def\vr{\varrho}
\def\vt{\theta}
\def\z{\zeta}

\def\grad{{{\mathop{\rm{grad}}}\,}}

\def\what{\widehat}

\def\cn{\bC^n}

\def\oD{\overline D}

\def\phim{\phi^-_1}

\def\phip{\phi^+_1}

\def\Z{\mathbb Z}

\def\R{\mathbb R}

\def\C{\mathbb C}

\def\tphi{(\theta,\phi)}

\def\tdt{\tilde\tau}

\def\rtor{\mathbb R\rightarrow\mathbb R}

\def\Rm{\begingroup\let\par=\null\obeylines\RemovE}
\def\RemovE#1\mR{\endgroup}

\begin{document}
\title{Hulls of Surfaces}
\author[A.\ J.\ Izzo]{Alexander J. Izzo}
\thanks{This paper was completed while the first author was a visitor at Indiana University.  He would like to thank the Department of Mathematics for its hospitality. }
\address{Alexander J. Izzo, Department of Mathematics and Statistics, Bowling Green State University, Bowling Green, OH, 43403 USA}
\email{aizzo@math.bgsu.edu}
\author[E.\ L.\ Stout]{Edgar Lee Stout}
\address{Edgar Lee Stout, Department of  Mathematics, University of Washington, Seattle, Washington, 98195.
}
\email{stout@math.washington.edu}

%
%
\date{\today}

\begin{abstract}
In this paper it is shown that every compact two-dimensional manifold $S$, with or without boundary, can be embedded in $\bC^3$ as a smooth submanifold $\Sigma$ in such a way that the polynomially convex hull of $\Sigma$, though strictly larger than $\Sigma$, contains no
analytic disc.

\end{abstract}

\maketitle

\section{ {Introduction.}}

The main result of the present note is the following theorem.

\begin{theorem}
\label{3.31.14.iii}
If  $S$ is a compact surface, i.e., a compact two-dimensional manifold, then there is a smooth surface $\Sigma$ embedded in $\bC^3$ that is homeomorphic to $S$ and that has the property that the polynomially convex hull $\what \Sigma$, although  strictly larger than $\Sigma$, contains no analytic disc.  Furthermore, the surface $\Sigma$ can be chosen to be totally real.
\end{theorem}

We note explicitly that the surfaces contemplated in the theorem may be closed, i.e., without boundary, or with boundary. The case of surfaces with boundary will be seen to be an immediate consequence of the case of closed surfaces.

Theorem 1.1 gives an example of a set, the surface $\Sigma$, with the property that its polynomially convex hull contains no analytic disc. Early in the study of uniform algebras and polynomial convexity it was conjectured that if the polynomially convex hull $\what E$ of a compact set $E$ in $\bC^n$ is strictly larger than $E$, then the complementary set $\what E\setminus E$ must contain an analytic disc. This optimistic expectation was shown to be wrong by Stolzenberg \cite{Stolzenberg:1963a} whose example is a suitable limit of analytic varieties. Since the appearance of Stolzenberg's counterexample, a variety of additional examples have been constructed. In \cite{Basener:1973} Basener constructed a smooth $3$--sphere in $\mathbb C^6$ that is not polynomially convex but whose hull has no analytic structure. Other examples of manifolds with this property are given in the paper \cite{I-SK-W:2020}. It is known from Stolzenberg's work on the hulls of curves \cite{Stolzenberg:1966} that for a smooth curve $K$ in $\mathbb C^n$, the set $\what K\setminus K$, if not empty,
is a purely one-dimensional variety.

In this paper we understand {\it{smooth}} functions or manifolds to be of class $\sC^\infty$. In addition, we will consider only connected surfaces.

The interest in the observation that the surface $\Sigma$ in the theorem can be taken to be totally real stems from a theorem of Alexander~\cite{Alexander:1996}: if $\Sigma$ is a totally real $n$--dimensional smooth submanifold of $\cn$, the complementary set $\what\Sigma\setminus\Sigma$, which is known a priori to be nonempty, necessarily contains an analytic disc.

Theorem \ref{3.31.14.iii} is an extension of a result in
 \cite{I-SK-W:2020} that exhibits a smooth two-dimensional torus $T^2$ in $\bC^3$  whose hull, which is strictly larger than $T^2$, contains no analytic disc.

By an {\it{analytic disc}} in $\cn$ we understand a set of the form $g(\bU)$ with $g$ a nonconstant holomorphic $\bC^n$--valued function on the open unit disc $\bU$ in $\bC$. Occasionally we speak of $g$ itself as being an analytic disc. An analytic disc $g(\bU)$ can be topologically complicated: For example, every irreducible one-dimensional analytic subvariety of an open set in $\cn$ is an analytic disc in our sense.
(That this is so is easily seen if not familiar: If $V$ is an irreducible one-dimensional subvariety of an open subset of $\cn$, let $(\mathcal R,\eta)$ be its normalization so that
$\mathcal R$
is
 a connected Riemann surface and $\eta:\mathcal R\rightarrow V$ is a surjective holomorphic map
that has certain additional properties.
The universal covering space   $(\mathcal R^*,
\eta^*)$
 of the surface $\mathcal R$ is the plane $\bC$ or the unit disc
$\bU$. In the latter case, the map $\eta\circ\eta^*:\bU\rightarrow V$ exhibits $V$ as an analytic disc in our 
sense; in the former case, let $f$ be a holomorphic function on the disc with range all of $\bC$. The map $\eta\circ\eta^*\circ f$ exhibits $V$ as an analytic disc.)

Theorem \ref{3.31.14.iii} is derived using the method of connected sums from the following result. We use the notation that $V_e$ denotes the totally real $2$--dimensional plane given by
\begin{equation}
\label{3.29.16.i}V_e=\{(x_1+ix_2, 2x_2-ix_1): x_1, x_2\in \bR\}
\end{equation}
and that $\bB_n(r)$ is the open ball in $\cn$ of radius $r$ centered at the origin. 

\begin{theorem}
\label{9.3.14.i} The space $\bC^2$ contains a smooth closed  submanifold $\Delta$ that is diffeomorphic to the open unit disc in the complex plane and
that has these properties:

\noindent(i) $\Delta$ contains a compact subset $E$ with the property that $\what E\setminus E$, though nonempty, contains no analytic disc.

\noindent(ii) There is  an $R>0$ such that the 
part of $\Delta$ outside the ball $\bB_3(R)$ coincides with the part of $V_e$ outside $\bB_3(R)$.

\noindent(iii) There is a smooth complex-valued function $f$ on $\Delta$ with the properties   that $f^{-1}(0)=E$ and that
each compact subset of each level set $f^{-1}(z)$ for $z\in\bC\setminus \{0\}$ is polynomially convex
and satisfies $\sP(f^{-1}(z))=\sC(f^{-1}(z))$.

\noindent(iv) The set $f^{-1}(1)$ is contained in $V_e$ and is the 
closure of the unbounded component of $\Delta\setminus C$ for a simple closed curve $C$ contained in  $\Delta \cap \bigl(\bC^2\setminus \bB_3(R)\bigr)$.

\noindent(v) The range of $f$ is contained in the subset $[-1,1]\cup i\bR$ of the coordinate axes in the complex plane.

\noindent(vi) The graph of $f$ is totally real. 
\end{theorem}

\smallskip
 
 Theorem~\ref{9.3.14.i} will be seen to be a consequence of a result of Alexander \cite{Alexander:1998} according to which the standard torus 
 $\bT^2=\{(e^{i\vt},e^{i\phi}):\vt,\phi\in\bR\}$ in $\bC^2$
 contains a compact subset $E$ such that $\what E\setminus E$ is not empty but contains no analytic disc. Moreover, if $V$ is any analytic subvariety of the bicylinder $\bU^2$ whose boundary $bV=\bar V\setminus V$ is contained in $\bT^2$ and $W$ is a neighborhood of $bV$ in $\bT^2$, then  such a set $E$ can be found in $W$. 
 In particular, there is such a set $E$ in every neighborhood of the diagonal in $\bT^2$.
\smallskip

\smallskip
Before beginning the proofs of the theorems above certain preliminaries are required.

\section{Preliminaries on Surfaces.}

Recall the classification of compact closed surfaces:

\begin{theorem}
\label{3.12.14.i}
A compact closed surface is of one of three forms: If it is orientable, it is a sphere or the connected sum of a finite number of tori. If it is not orientable, it is the connected sum of a finite number of real projective planes.
\end{theorem}

Given surfaces $S_1$ and $S_2$, their connected sum is denoted by $S_1\#S_2$ and is defined to be the surface obtained by excising an open disc with smooth boundary from each of $S_1$ and $S_2$ to obtain surfaces $S'_1$ and $S'_2$  and gluing the boundaries of $S'_1$ and $S'_2$, each of which is a smooth simple closed curve, together with a diffeomorphism. In the event that both $S_1$ and $S_2$ are compact, the result is another compact surface.

The binary operation $\#$ is commutative and associative in the sense that $S_1\#S_2$ is homeomorphic to $S_2\#S_1$ and, given a third surface $S_3$, the surface $S_1\#(S_2\#S_3)$ is homeomorphic to the surface $(S_1\#S_2)\#S_3$. If $S$ is any surface, then the connected sum of $S$ and a $2$--sphere is homeomorphic to $S$.

The theory of compact surfaces is developed in detail in the books of Massey \cite{Massey:1967} and of Munkres \cite{Munkres:2000}.

\section{Preliminaries on Polynomial Convexity.}

Recall that if $X$ is a compact subset of $\cn$, its polynomially convex hull $\what X$ is the compact subset of $\cn$ defined by
$$\what X=\{z\in\cn:|P(z)|\leq \max_{x\in X}|P(x)|\ 
\mbox{\rm{for\, all\, polynomials}}\, P
\}.$$

For a compact set $X$ in $\cn$, we denote by $\sP(X)$ the uniformly closed subalgebra of $\sC(X)$ generated by the holomorphic polynomials.

In general, given $X$ it is difficult to determine $\what X$. There is, though, the following criterion, which will be useful in what follows.
\begin{lemma}{\rm {\cite[pp.13\,-14.]{Stout:2007}}}
\label{3.30.14.iv}
If $X$ is a compact subset of $\cn$ and if $g\in\sP(X)$ is such that $\sP(g(X))=\sC(g(X))$, then $X$ is polynomially convex if and only if each of the fibers $X_t=g^{-1}(t)$, $t\in\bC$, is polynomially convex. If $X$ is polynomially convex, then $\sP(X)=\sC(X)$ if and only if for each $t\in\bC$, the fiber $X_t$ satisfies $\sP(X_t)=\sC(X_t)$.
\end{lemma}

This lemma has the following immediate corollary:

\begin{corollary}
\label{3.30.14.vi}
Let $X$ be a compact subset of $\cn$, and let $X_0$ be a closed subset of $X$. If $f\in\sC(X)$ is a complex-valued function with $X_0=f^{-1}(0)$ that satisfies the conditions that $\sP(f(X))=\sC(f(X))$ and  that each fiber $X_t=f^{-1}(t),\, t\in\bC\setminus\{0\}$, is polynomially convex, then the graph $\Gamma_f$ of $f$ satisfies $\what\Gamma_f=(\what X_0\times \{0\})\cup\Gamma_f\subset \bC^{n+1}$. 
\end{corollary}

There is a further result in the spirit of the preceding results. 

\begin{lemma}
\label{6.18.14.ii}
Let $X$ be a compact subset of $\cn$, and let $X_0$ be a closed subset of $X$. If $f\in\sC(X)$ is a complex-valued function with $X_0=f^{-1}(0)$ that satisfies the conditions that $\sP(f(X))=\sC(f(X))$ and that $\sP(f^{-1}(t))=\sC(f^{-1}(t))$ for all $t\in\bC\setminus\{0\}$, then $\sP(\Gamma_f)=\{g\in\sC(\Gamma_f):g|(X_0\times\{0\})\in\sP(X_0\times\{0\})\}$.
\end{lemma}
\noindent{\bf Proof.} If $E$ is a compact subset of 
$f(X)\setminus\{0\}$ then the set 
$$E_f=\{(z,f(z)):z\in X\ \mbox{\rm{and}}\ f(z)\in E\}$$
 is a compact subset of $\Gamma_f$ that is a peak set for $\sP(\Gamma_f)$ and that satisfies $\sP(E_f)=\sC(E_f)$ as follows from Lemma~\ref{3.30.14.iv}. Accordingly, if $\mu$ is a measure on $X$ that annihilates the algebra $\sP(X)$, then $\mu$ is concentrated on $X_0$, whence the lemma.

We will need below the following result, which is Lemma 3.2 from the paper \cite{I-SK-W:2020}.
\begin{lemma}
\label{2.26.15.i}
 Let $K$ be a closed subset of $\bT^2$  such that for some complex number $c$ of modulus one the
set $K$ is disjoint from the circle $\{ (z_1,z_2)\in\bT^2 : z_1=c\}$ and there is no
$a$ of modulus one such that $K$ contains the full circle $\{(z_1, z_2)\in\bT^2 : z_1 = a\}$. Then
$\sP(K) = \sC(K)$, and, in particular, $K$ is polynomially convex.
\end{lemma}

\section{Preliminaries on Totally Real Embeddings.}

In the proof of Theorems \ref{3.31.14.iii} and   \ref{9.3.14.i} we need some information about totally real embeddings.

Given a smooth surface $\Sigma$ in $\bC^2$ and a smooth function $f$ on $\Sigma$, we ask when the graph $\Gamma_f=\{(z_1,z_2,f(z_1,z_2)):(z_1,z_2)\in\Sigma\}$ is totally real.
 An answer can be given as follows: Let $\Phi:\Sigma\rightarrow\bC^3$ be the graph map given by 
\begin{equation}
\label{6.20.14.ii}
\Phi(z_1,z_2)=(z_1,z_2,f(z_1,z_2)).
\end{equation}
\begin{lemma}
\label{6.20.14.i}
The graph $\Gamma_f $ is totally real if and only if for every point $p\in\Sigma$, there is a holomorphic two-form $\varTheta_p$ on $\bC^3$ whose pull-back $\Phi^*\varTheta_p$, a smooth two-form on $\Sigma$, does not vanish at $p$.
\end{lemma}
Alternatively put, $\Gamma_f$ is totally real at a point $\Phi(p)$  if and only if one of the forms $dz_1\wedge dz_2,\ dz_1\wedge df,\ dz_2\wedge df$ does not vanish at $p$.

As a consequence of this lemma, we see that if we write $\Sigma=\Sigma_\bC\cup\Sigma_\bR$ with $\Sigma_\bC$
the closed set of points in $\Sigma$ at which the two-plane tangent to $\Sigma$ is a complex line
and with $\Sigma_\bR$ the open subset of $\Sigma$ of points at which the tangent plane is a totally real two-plane, then $\Gamma_f$ is totally real at every point of $\Phi(\Sigma_\bR)$: With holomorphic coordinates $w_1,w_2,w_3$ on $\bC^3$, $\Phi^*(dw_1\wedge dw_2)=dz_1\wedge dz_2$ vanishes at no point of $\Sigma_\bR$.
The problem of the total reality of $\Gamma_f$ is thus seen to be that of ascertaining whether the
graph of $f$ is totally real over the points of $\Sigma_\bC$.

 In particular, if $\Sigma$ is totally real, then no matter what the smooth function $f$ may be, the graph $\Gamma_f$ is a totally real submanifold of $\bC^3$. 
For example, the standard torus in $\bC^2$ is totally real, so the graph of every smooth function on it is again totally real.

For graphs of functions on surfaces not assumed to be totally real, there is a simple sufficient condition: If $f$ is real-valued and smooth on the surface $\Sigma$ in $\bC^2$ and if, moreover, $df$ does not vanish, then the graph $\Gamma_f$ is totally real over any point of $\Sigma$ at which the tangent is a complex line.

We now turn to the second of the theorems of the Introduction.

\section{Proof of Theorem~\ref{9.3.14.i}.}

 Fix $0<\alpha_1<\alpha_2<\pi/2<\alpha_3<\alpha<\pi$.  Choose $\sC^\infty$-smooth  real-valued functions $\beta$ and $\gamma$ on $\bR$ and $\chi$ on $(-\alpha,\infty)$ such that 
\begin{equation*}
\begin{split}
\beta(\phi)&=\begin{cases}&0\ \mbox{for}\ |\phi|\leq\alpha_1\\
&1\ \mbox{for}\  |\phi|\geq \alpha_3\\
\end{cases}\\
\gamma(\phi)&=\begin{cases}&\phi\ \mbox{for}\ |\phi|<\alpha_2\\
&\pi/2\ \mbox{for}\  \phi\geq \alpha_3\\
&-\pi/2\ \mbox{for}\ \phi\leq -\alpha_3\\
\end{cases}\\
\chi(\phi)&=\begin{cases}&1\ \mbox{for}\ |\phi|<\alpha_2\\
&0\ \mbox{for}\  \phi\geq \alpha\\
\end{cases}\\
\end{split}
\end{equation*} We require also that $\lim_{\phi\rightarrow -\alpha^+}\chi(\phi)=\infty$, that $\beta$ be an even function, that $\gamma$ be an odd function, and that each of $\beta$, $\gamma$, and $\chi$ be strictly monotonic
 on each interval where the values have not been specified by the formulas above and, indeed, that the derivatives of these functions not vanish at any point of any of these intervals.

Let $\mathcal R$ be the half-open rectangle in the $(\vt,\phi)$--plane  given by
$$\mathcal R=[-\pi,\pi]\times (-\alpha,\alpha],$$ 
define $\Phi:\bR\times(-\alpha,\alpha]\rightarrow \bC^2$ by
\begin{equation}
\label{1.12.15.i}
\Phi(\vt,\phi)=\Bigl(\chi(\phi)e^{i\vt}, \chi(\phi)\bigl[e^{i(\gamma(\phi)+\vt)}+\beta(\phi)\cos\bigl(\gamma(\phi)+\vt\bigr)\bigr]\Bigr),
\end{equation}
and set
$$\Delta=\Phi(\mathcal R).$$
The map  $\Phi$ 
 is proper, and one verifies that it 
is injective on $\mathcal R$ except for identifying each point $(-\pi, \phi)$ with the point $(\pi,\phi)$ and shrinking the interval $[-\pi,\pi]\times\{\alpha\}$  
to the point $(0,0)$.  Thus $\Delta$ is topologically an open disc.

We claim that, furthermore, the map $\Phi$ 
 is an immersion on the set $\bR\times(-\alpha,\alpha)$, that is, has injective derivative at each point.  To verify this we work with real coordinates.  If $x_1,\ldots, x_4$ are the real coordinate functions of the map $\Phi$, then,
$$\quad\frac{\partial(x_1,x_2)}{\partial(\vt,\phi)}=
\left[   \begin{array}{rr}
-\chi(\phi)\sin\vt & \chi'(\phi)\cos\vt \\*[4pt]
\chi(\phi) \cos\vt & \chi'(\phi) \sin\vt
\end{array} \right]$$
so that $\det\left(\displaystyle{\frac{\partial(x_1,x_2)}{\partial(\vt,\phi)}}\right)=\,-\chi(\phi)\chi'(\phi)\neq 0$ \vadjust{\kern 6pt} whenever $\alpha_2<|\phi|<\alpha$; consequently, 
$\displaystyle{\frac{\partial(x_1,x_2,x_3,x_4)}{\partial(\vt,\phi)}}$ has rank two for these $\phi$.
For $|\phi|\leq \alpha_2$, note that 
$$x_3+ix_4=e^{i(\phi+\vt)}+\beta(\phi)\cos(\phi+\vt)=
\bigl(1+\beta(\phi)\bigr)\cos(\phi+\vt) +i\sin(\phi+\vt),$$ whence
$$\frac{\partial(x_3,x_4)}{\partial\phi}=
\left[   \begin{array}{c}
-\bigl(1+\beta(\phi)\bigr)\sin(\phi+\vt) + \beta'(\phi)\cos(\phi+\vt)  \\*[4pt]
\cos(\phi+\vt)
\end{array} \right].$$
Since $\cos(\phi+\vt)$ and $\sin(\phi+\vt)$ can not both be simultaneously zero, and $1+\beta(\phi)$ is never zero, we conclude that $\displaystyle\frac{\partial(x_3,x_4)}{\partial\phi}$ is \vadjust{\kern 6pt} never zero.  Since $\chi(\phi)=1$ for all $|\phi|\leq \alpha_2$, we get that for all $\phi$ in that \vadjust{\kern 6pt}interval
$$\quad\frac{\partial(x_1,x_2, x_3,x_4)}{\partial(\vt,\phi)}=
\left[   \begin{array}{cc}
\sin\vt & \quad0 \\*[4pt]
\cos\vt & \quad0 \\*[4pt]
* & \quad A \\*[4pt]
* & \quad B \\*[4pt]
\end{array} \right]$$
where $(A,B)\neq (0,0)$.  Thus $\displaystyle{\frac{\partial(x_1,x_2,x_3,x_4)}{\partial(\vt,\phi)}}$ has rank two for these $\phi$ as well.

From the result of the preceding paragraph, it follows that $\Delta$ is smooth everywhere except possibly at the origin.  However, $\Delta$ contains a neighborhood of the origin in the (totally real) plane $$V_i=\{(x_1+ix_2, -2x_2+ix_1): x_1, x_2\in \R\}$$ and hence is smooth at the origin also:  If $\mathcal R_{\alpha_3}^+$ is the closed rectangle in the $(\vt,\phi)$--plane given by
$$\mathcal R_{\alpha_3}^+=[-\pi,\pi]\times[\alpha_3,\alpha],$$
then $\Phi(\mathcal R^+_{\alpha_3})\subset V_i.$

For $\tau\in (0,\pi)$, let 
$A_{\tau}$ denote the topological annulus on the torus $\bT^2$ given by
$$A_{\tau}=\{(e^{i\vt}, e^{i(\phi + \vt)}): -\pi\leq\vt\leq \pi, -\tau< \phi<\tau\}.$$ With this notation, 
$A_{\alpha_1}$ is contained in $\Delta$ and is a neighborhood of the diagonal in $\bT^2$.

Outside a Euclidean $3$-sphere of some sufficiently large radius $R$ centered at the origin of $\bC^2$, the disc $\Delta$ is the unbounded annulus consisting of points with norm bigger than $R$ lying in the totally real plane $V_e$ defined in equation (\ref{3.29.16.i}).

In fact, if 
\begin{equation}
\label{3.29.16.iii}
\mathcal R_{-\alpha_3}^-=[-\pi,\pi]\times(-\alpha,-\alpha_3]
\end{equation}
then $\Phi(\mathcal R_{-\alpha_3}^-)$ is an unbounded annulus contained in $ V_e$.
\smallskip

We now construct the function $f$. The construction depends on Lemma~3.5 of \cite{I-SK-W:2020}.

Choose angles $\phim$, $\phip$, $\phi_2$ with $0<\phim<\phip<\alpha_1<\phi_2<\alpha_2$, and with $\phip-\phim=2\pi/n$ for some positive integer $n$ divisible by 4.  Let $E_1$ be a compact set contained in the annulus $\{(e^{i\vt}, e^{i(\phi + \vt)}): -\pi\leq\vt\leq \pi, -\phim< \phi<\phim\}$ on $\mathbb T^2$ such that $\what E_1\neq E_1$ but $\what E_1$ contains no analytic disc.
The existence of such a set $ E_1$ is provided by \cite{Alexander:1998} as noted above at the end of the introduction.

Let $Z_+$ be the sawtooth path in $\R^2$ consisting of the vertical line segments 
$$\{(2\pi m/n, \phi): m\in \Z, \phim\leq\phi\leq\phip\}$$ together with the diagonal line segments connecting, for each $m$, the pair of points \hbox{$(2\pi m/n, \phip)$} and $(2\pi (m+1)/n, \phim)$.  Let $Z_-$ be the sawtooth path in $\R^2$ that is the image of $Z_+$ under the map $(\vt,\phi)\mapsto (-\vt, -\phi)$.
Let $L$ consist of two vertical line segments defined by
$$L=\{0,\pi\}\times [-\phip,\phip].$$
Finally define a subset $E$ of $\bT^2$ by
$$E=E_1\cup \Phi(Z_+\cup Z_-\cup L)$$
with $\Phi$ as given in equation~(\ref{1.12.15.i}).
Note that $\Phi(Z_+\cup Z_-\cup L)$ is the union of a finite set of circular arcs on $\bT^2$ each of which has a constant $z_1$ or $z_2$ value.  Thus repeated application of Lemma~3.5 of \cite{I-SK-W:2020} yields 
$$\what E\setminus E=\what E_1\setminus E_1$$
and hence $\what E\neq E$ but $\what E$ contains no analytic disc.

Fix $\phi_3$ with $\alpha_3<\phi_3<\alpha$.    
Let $G_o$ denote the closed set of points in $\R^2$ that lie on or below the sawtooth path $Z_+$, and let $G_u$ denote the open set of points in $\R^2$ that lie strictly above the sawtooth path $Z_+$ and strictly below the line $\phi=\phi_3$.  Let $W=\{(\vt,\phi)\in \R^2: \frac{\pi}{2}-\frac{\pi}{2n}<\vt<\frac{\pi}{2}+\frac{\pi}{2n},\  \phip<\phi<\phi_2\}$.

\begin{lemma}\label{2.15.15.i} There exist $\sC^\infty$ functions $g,h: \R^2\rightarrow \R$ such that
\hfil\break
(i) $g\tphi=0$ for $\tphi\in G_o$\hfil\break
(ii) $(\partial g/\partial \phi)\tphi >0$ for $\tphi\in G_u$\hfil\break
(iii) $g\tphi$ is independent of $\vt$ for $\phi\geq \phi_2$\hfil\break
(iv) $g\tphi=1$ for $\phi\geq \phi_3$ (and all $\vt$)\hfil\break
(v) $g\tphi$ is $2\pi/n$ periodic in $\vt$ for fixed $\phi$\hfil\break
(vi) conditions (i)--(iv) hold also with  $g$ replaced by $h$\hfil\break
(vii) $h\tphi=g\tphi$ for all $\tphi\in \R^2\setminus W$\hfil\break
(viii) $h\tphi>g\tphi$ for all $\tphi\in W$.\hfil\break
\end{lemma}
Note that condition (v) fails with $g$ replaced by $h$.

\smallskip

Assume the Lemma for the moment.  The functions $g$ and $h$ induce $\sC^\infty$ functions $\tilde g$ and $\tilde h$ on $\Delta$ via the parametrization $\Phi$  restricted to $[-\pi,\pi]\times (-\alpha,\alpha]$.  Define $h_r$ on $\R^2$ by $h_r\tphi=h(-\vt,-\phi)$, and let $\tilde h_r$ be the induced function on $\Delta$.  Let $H_u$ denote $G_u \cup \{\tphi: \phi_3\leq\phi\leq\alpha\}$ so that $H_u$ consists of everything above the sawtooth  path $Z_+$ up to $\phi=\alpha$. Also, let $H'_u$ be the subset of $H_u$ comprising those points $(\vt,\phi)$ with $\phi<\alpha$.

  Choose, by Lemma~3.3 of \cite{I-SK-W:2020}, a smooth real-valued function $k$ on $\Delta$  that is identically zero on the closure of the set $E\cup \Phi(H_u)\cup \Phi(-H'_u)$,
is strictly positive on the part of the complement of this set corresponding to $\vt$s satisfying  
 $0<\vt<\pi$, and is strictly negative on the part of the complement of this set  corresponding to $\vt$s satisfying
 $-\pi<\vt<0$.  Then set
$$f=\, -\tilde h + \tilde h_r + ik.$$
 
\smallskip

The function $f$ is constantly one at the points $\Phi(\vt,\phi)$ with $\phi\in(-\alpha,-\phi_3)$, so for the simple closed curve $C$ of the statement of the Theorem, we can take the  set $\{\Phi(\vt,-\phi_3):\vt\in[-\pi,\pi]\}$.

The range of the function $f$ is contained in the union $\bR\cup i\bR$ of the coordinate axes whence for every compact subset $X$ of $\Delta$, $\sP(f(X))=\sC(f(X))$. Note that $f^{-1}(0)=E$, and recall that 
we have already observed that $\what E\setminus E$ is not empty but has no analytic structure.  Each of the level sets where $f$ is pure imaginary is polynomially convex by Lemma~3.2 of \cite{I-SK-W:2020}, which was quoted above as Lemma~\ref{2.26.15.i}.  Each of the level sets of $f$ in the regions corresponding to values of $\phi$ with 
$\phi\geq \phi_2$ or $\phi\leq \, -\phi_2$ is contained in a totally real plane and hence is polynomially convex. Each of the level sets of $f$ corresponding to values of $\phi$ in the range
 $-\phi_2<\phi<\phi_2$ on which 
$f$ is non-zero and real is polynomially convex because it is the graph of a function on the unit circle that does not belong to the disc algebra (and hence is polynomially convex by Wermer's maximality theorem). To establish this last assertion we argue as follows.

Consider a level set of $f$ where $f$ is non-zero and real that corresponds to a $\phi$ with
 $0<\phi<\phi_2$.  (The argument for $-\phi_2<\phi<0$ is the same.)  This set is $\{f=\, -c\}$
for some $c>0$.  Note that $\{f=\, -c\}$  is $\Phi(\{h=c\})$.  Also 
$$\{h=c\}\setminus W\,=\,\{g=c\}\setminus W.$$
Because $(\partial g/\partial \phi)\tphi >0$ for $\tphi\in G_u$, the set $\{g=c\}$ is a graph $\phi=u(\vt)$ for some smooth function $u$.  The same applies to the set $\{h=c\}$; it is a graph, say $\phi=v(\vt)$.  Also $v=u$ outside of the interval $\frac{\pi}{2}-\frac{\pi}{2n}<\vt<\frac{\pi}{2}+\frac{\pi}{2n}$ because $h=g$ outside of $W$.  Because $g$ is $2\pi/n$ periodic in $\vt$, so is $u$.  The set $\{f=\, -c\}=\Phi(\{h=c\})$ is the graph over the unit circle of the function $V$ defined by
$$V(e^{i\vt})= e^{i(\gamma(v(\vt))+\vt)}+\beta(v(\vt))\cos\bigl(\gamma(v(\vt))+\vt\bigr).$$
Let $U$ denote the function defined by the same formula as $V$ but with $v$ replaced by $u$.
Because $v=u$ outside the interval $\frac{\pi}{2}-\frac{\pi}{2n}<\vt<\frac{\pi}{2}+\frac{\pi}{2n}$, we get that for $\theta$ outside that interval $V(e^{i\vt})=U(e^{i\vt})$.  A  computation shows that the $2\pi/n$ periodicity of $u$ implies that 
\begin{equation}\label{2.19.15.i}
U(e^{i(\vt+\pi)})= \, - U(e^{i\vt}).
\end{equation}
Now assume to get a contradiction that the set $\{f=\, -c\}=\Phi(\{h=c\})$ is the graph of a disc algebra function, i.e., that $V$ is in the disc algebra.  Then $U$ agrees with a disc algebra function on the circular arc $\{e^{i\vt}:-\pi<\vt< \pi/4\}$.  Equation~(\ref{2.19.15.i}) then implies that $U$ also agrees with a disc algebra function on the circular arc $\{e^{i\vt}:0<\vt <5\pi/4)$.  Since these arcs overlap and cover the unit circle, it follows that $U$ must itself be in the disc algebra.  But then $V$ cannot be in the disc algebra since $V$ agrees with $U$ on an arc but is distinct from $U$.  This contradicts our assumption and hence establishes the polynomial convexity of the level set $\{f=\, -c\}$.

We now know that compact subsets of  all the level sets of $f$ save only $f^{-1}(0)$ are polynomially convex. Moreover, the discussion above  makes it evident that each compact subset $K$ of a level 
$f^{-1}(z),\ z\neq 0$, satisfies $\sP(K)=\sC(K)$.  
\medskip

We finally verify that the graph of the function $f$ is totally real. To this end, recall that the subset $\Phi(R^-_{-\alpha_3})$ is contained in the totally real plane 
$V_e$, and the subset 
$\Phi(R^+_{\alpha_3})$ is contained in the totally real plane ${V_i}$.  Moreover, the part of $\Delta$ comprising the points $\Phi(\theta,\varphi)$ with 
$|\varphi| <1$ lies in the totally real torus $\bT^2$.  Thus, the part of the graph of 
$f$ that lies over these sets is totally real. (Recall the remarks immediately following Lemma \ref{6.20.14.i}.)  For the part of the graph of $f$ that lies over the remainder of $\Delta$, note that $f$ is real-valued at all points 
$\Phi(\theta, \varphi)$ with $\varphi_1^+ < \varphi<\varphi_3$ and that $df$ is nonzero there since $\partial (f \circ \Phi)/\partial \varphi$ is strictly positive through this region.  Consequently, this part of the graph of $f$ is totally real as well---again recall the remarks following Lemma \ref{6.20.14.i}.

The function $f$ has now been shown to have the required properties, and the proof of Theorem~\ref{9.3.14.i} is complete except for the proof of the lemma.

\smallskip

\noindent{\bf Proof of Lemma~\ref{2.15.15.i}}.  Let $K_u=G_u\cup \{\tphi:\phi\geq\phi_3\}$, i.e., $K_u$ consists of everything above the sawtooth path.  Set $W_1=\{(\vt,\phi)\in \R^2: \frac{\pi}{2}-\frac{\pi}{2n}<\vt<\frac{\pi}{2}+\frac{\pi}{2n},\  \phi>\phip\}$ so that $W=W_1\cap \{\phi<\phi_2\}$, and choose a countable collection $W_2, W_3, \ldots$ of open sets of the form $W_i=(r_i,s_i)\times (t_i,\infty)$ with $r_i,s_i,t_i\in \R$ such that $W_i\subset K_u$ for every $i$ and $\{W_i\}_{i=1}^\infty$ covers the set $K_u\cap \bigl([\frac{\pi}{2}-\frac{\pi}{n}, \frac{\pi}{2}+\frac{\pi}{n}]\times \R\bigr)$. 
We require also that $s_i-r_i\leq\pi/n$ for all $i$.  For each $j\in \Z$ set  
$$W_{ij}=(2\pi j/n, 0)+W_i,$$
which is the translate of $W_i$  
by $2\pi j/n$ in the $\vt$-direction.  Then
$$\bigcup_{ij}W_{ij}=K_u.$$

For each $i=1,2,\ldots$, choose a $\sC^\infty$ function $\sigma_i:\bR\rightarrow\bR$ with $\sigma_i>0$ on the interval $(r_i,s_i)$ and 
$\sigma_i=0$ outside $(r_i,s_i)$.
Choose a $\sC^\infty$ function $\tilde\tau_0:\bR\rightarrow\bR$ such that $\tdt_0(\phi)=0$ for $\phi\leq\phip$, $\tilde\tau_0(\phi)=1$ for $\phi\geq\phi_2$, and $\tilde\tau'_0(\phi)>0$ for $\phip<\phi<\phi_2$.  Set $\tilde\tau_1=\tilde\tau_0^2$.  Note that $\tilde\tau_1=\tilde\tau_0$ on $(-\infty, \phip]\cup[\phi_2,\infty)$ and $\tilde\tau_0>\tdt_1$ on $(\phip,\phi_2)$.  For each $i=2,3,\ldots$, choose a $\sC^\infty$ function $\tdt_i:\bR\rightarrow\bR$ such that 
$\tdt_i(\phi)=0$ for $\phi\leq t_i$, $\tdt_i(\phi)=1$ for $\phi\geq\phi_2$, and 
$\tdt'_i(\phi)>0$ for $t_i<\phi<\phi_2$.  Finally choose a $\sC^\infty$ function $\tau:\rtor$ such that $\tau(\phi)=0$ for $\phi\leq 0$, $\tau(\phi)=1$ for $\phi\geq \phi_3$, and $\tau'(\phi)>0$ for $0<\phi<\phi_3$.

For each $j\in \Z$, define $\sigma_{ij}:\rtor$ by
$$\sigma_{ij}(\vt)=\sigma_i(\vt-2\pi j/n),$$
i.e., $\sigma_{ij}$ is the translate of $\sigma_i$ by $2\pi j/n$.  Then at every point  $\vt\in \R$ there is some $\sigma_{ij}$ that is strictly positive there.  Also for fixed index $i$, the support of $\sigma_{ij}$ is disjoint from that of $\sigma_{ij'}$ for $j\neq j'$.

For each $i=0,1,2,\ldots$, set
$$\tau_i=\tdt_i\tau.$$
For notational convenience, set $t_0=t_1=\phip$.
Note that then for each $i=0,1,2,\ldots$, \hfil\break
(i$^\prime$) $\tau_i(\phi)=0$ for $\phi\leq t_i$\hfil\break
(ii$^\prime$) $\tau'_i(\phi)>0$ for $t_i<\phi<\phi_3$\hfil\break
(iii$^\prime$) $\tau_i(\phi)=\tau(\phi)$ for $\phi\geq\phi_2$ \hfil\break
(iv$^\prime$) $\tau_i(\phi)=1$ for $\phi\geq\phi_3$.\hfil\break
Note also that $\tau_1=\tau_0$ on $(-\infty, \phip]\cup[\phi_2,\infty)$ and $\tau_0>\tau_1$ on $(\phip,\phi_2)$.

\smallskip

If the sequence $\{c_i\}_{i=1}^\infty$ of positive numbers decreases sufficiently rapidly as $i\rightarrow\infty$, then the series

$$\tilde g\tphi=\sum_{i=1}^\infty\sum_{j=-\infty}^\infty c_i\sigma_{ij}(\vt)\tau_i(\phi)$$
converges and its sum, which we denote by $\tilde g(\vt,\phi)$, is a $\sC^\infty$ function on $\bR^2$. As concerns the convergence of this series,
notice that the summands of the sum $\sum_{j=-\infty}^\infty \sigma_{ij}(\theta)$ have mutually disjoint support
and that, granted that $c_i$ decreases sufficiently rapidly, the sum over $i$, along with all its derivatives, converges uniformly.
For $\phi\geq\phi_2$
$$\tilde g\tphi=
\tau(\phi)\sum_{i=1}^\infty\sum_{j=-\infty}^\infty c_i \sigma_{ij}(\vt)$$
by condition (iii$^\prime$).  Define $\lambda:\rtor$ by
$$\lambda(\vt)=\sum_{i=1}^\infty\sum_{j=-\infty}^\infty c_i\sigma_{ij}(\vt).$$
Note that $\lambda(\vt)>0$ for all $\vt$. Finally define $g:\R^2\rightarrow \R$ by
$$g\tphi=\tilde g\tphi/\lambda(\vt).$$
Then $g$ has all the required properties.

\smallskip

For $h$, define  
$\tilde h:\R^2\rightarrow \R$ by
$$\tilde h=c_1\Bigl(\sigma_{10}(\vt)\tau_0(\phi)+\sum_{j\neq 0}\sigma_{1j}(\vt)\tau_1(\phi)\Bigr)
+\sum_{i=2}^\infty\sum_{j=-\infty}^{\infty} c_i\sigma_{ij}(\vt)\tau_i(\phi).$$
(The function $\tilde h$ is given by the same double sum as $\tilde g$ except that the term $c_1\sigma_{10}(\vt)\tau_1(\phi)$ has been replaced by $c_1\sigma_{10}(\vt)\tau_0(\phi)$.)
Finally set
$$h\tphi=\tilde h\tphi/\lambda(\vt).$$
Then condition (vi) holds for the same reasons that conditions (i)--(iv) hold.  Conditions (vii) and (viii) hold on account of the relation between $\tau_0$ and $\tau_1$.

\section{{{The Proof of Theorem~\ref{3.31.14.iii}. First step.}}}

In this section we prove Theorem~\ref{3.31.14.iii} for a special class of surfaces, which eventually will be seen to include a topological copy of every closed surface.

\begin{lemma}
\label{2.27.15.i}
Let $\Sigma\subset \bC^2$ be a smooth two-dimenional compact submanifold without boundary with the following properties:

\noindent (i) There is a point $p\in\Sigma$ with a neighborhood $U$ in $\Sigma$ that is contained in the tangent plane ${\mbox{\bf T}}_p\Sigma$ of $\Sigma$ at the point $p$. 

\noindent (ii) The tangent plane $\mbox{\bf T}_p\Sigma$ is totally real.

\noindent (iii) There is a smooth real-valued function $g$ on $\Sigma$ with the property that $g=1$ on $U$ and 
 $g\geq 1$ on $\Sigma\setminus U$.

\noindent (iv) Each level set $g^{-1}(\z)$ for $\z\in\bR$ is polynomially convex.

Then there is a smooth submanifold $\Sigma^*$ of $\bC^3$ that is diffeomorphic to $\Sigma$ and that has the property that $\what{\Sigma^*}\setminus\Sigma^*$, although not empty, contains no analytic disc. 
If, in addition, $dg$ is nonzero at every point where $\Sigma$ has a complex tangent, then $\Sigma^*$ can be taken to be totally real.
\end{lemma}

\noindent{\bf Proof.}  Let  $\Delta$ be the disc of Theorem~\ref{9.3.14.i}, and let $f$ be the function on $\Delta$ given in Theorem~\ref{9.3.14.i}.

 We may suppose the point $p$ to be the origin and, by applying a nonsingular complex-linear (though not necessarily  unitary) automorphism of $\bC^2$, that $\mbox{\bf T}_p\Sigma$ is the totally real plane $V_e$ given by equation (\ref{3.29.16.i}).

Let $K\subset V_e$ be a circle that is centered at the origin. Require also that $K$ and the disc $D$ that it bounds in $V_e$ be contained in $U$.

It follows from conditions (ii) and (iv) of Theorem \ref{9.3.14.i} that if the positive number $\rho$ is sufficiently large, the expanded surface $\rho\Sigma$ will have the property that 
the unbounded component of $\Delta\setminus \rho K$ is contained in the set $f^{-1}(1)$.

 We now form a surface $\Sigma^{**}$ as follows: $\Sigma^{**}$ is the union of $(\rho\Sigma)\setminus\rho D$ and the part of $\Delta$ that lies inside the circle $\rho K$, i.e., the bounded component of $\Delta\setminus\rho K$. The surface $\Sigma^{**}$ so defined is smooth and is diffeomorphic to the initial surface $\Sigma$. 

Define the function $\tilde g$ on $\rho\Sigma$ by $\tilde g(z)=g(z/\rho)$.
Then  define a function $h$ on $\Sigma^{**}$ by the condition that on $(\rho\Sigma)\setminus\rho D$, the function $h$ coincide with the function $\tilde g$, and on the bounded component of $\Delta\setminus \rho K$, it coincide with the the function $f$ of Theorem~\ref{9.3.14.i}. This function is smooth, and its level sets $h^{-1}(\z)$ are polynomially convex for all $\z\neq 0
$. The fiber $h^{-1}(0)$ is the set $E$ of Theorem~\ref{9.3.14.i}.

It follows from Corollary \ref{3.30.14.vi} that the graph $\Sigma^*$ of the function $h$ is the surface we seek.  
The assertion about total reality of $\Sigma^*$ follows from the remark at the end of Section~4.
\smallskip

We now turn our attention to proving that every closed surface is homeomorphic to a surface in $\bC^2$ of the form described in the hypotheses of the preceding lemma. (Note that as every closed surface is homeomorphic to a smooth surface, we need only deal with smooth surfaces.)

\section{{Flattening Surfaces and Functions.}}

This section is devoted to a lemma about flattening surfaces and one about flattening functions.

We use consistently the notation that $\bB_n(c,r)$ denotes the open ball in $\bC^n$ centered at the point $c$ and of radius $r$ so that $\bB_n(r)=\bB_n(0,r)$.

\begin{lemma}
\label{4.16.14.i}
 Let $\Sigma$ be a compact smooth $n$-dimensional submanifold of $\bC^n$, let $p\in \Sigma$, and let $\Sigma$ be totally real near $p$.  There is an $\eta>0$ such that for each sufficiently small $\varepsilon>0$ there is a compact smooth submanifold $\Sigma_0$ of $\bC^n$ such that 
$\Sigma_0\setminus \bB_n(p,2\varepsilon)=\Sigma\setminus \bB_n(p,2\varepsilon)$, such that $\Sigma_0\cap \bB_n(p,\varepsilon)=\mbox{\bf T}_p(\Sigma)\cap \bB_n(p,\varepsilon)$ in which $\mbox{\bf T}_p(\Sigma)$ is the (totally real) $n$-plane tangent at $p$ to $\Sigma_0$, and such that $\sP(\Sigma_0\cap \overline{\bB_n(p,\eta)})=\sC(\Sigma_0\cap \overline{\bB_n(p,\eta)})$.  If the initial manifold $\Sigma$ is totally real, $\Sigma_0$ can be constructed to be totally real.
\end{lemma}

The proof will show, in addition, that if the initial submanifold $\Sigma$ lies in one of the closed half-spaces determined by a real hyperplane in $\bC^n$ that contains $\mbox{\bf T}_p(\Sigma)$, then the flattened surface $\Sigma_0$ lies in the same half-space.

\noindent{\bf Proof.} Let $\chi$ be a smooth nondecreasing function on the real line with $\chi(t)=0$ for $t\leq 1$ and $\chi(t)=1$ for $t\geq 2$. Let $\mu=\max_{0<t<2}\chi'(t)$. Define $\chi_\delta$ by $\chi_\delta(t)=\chi(t/\delta)$ so that $\chi_\delta=0$ for $t\leq \delta$ and $\chi_\delta=1$ for $t\geq 2\delta$. This function satisfies $\chi'_\delta(t)=\chi'(t/\delta)/\delta$ so that
\begin{equation}
\label{5.2.14.ii}
\chi_\delta'(t)\leq \mu/\delta\quad {\mbox{\rm whence}}\quad \chi_\delta(t)\leq \mu t/\delta\ {\mbox{\rm for}}\ t\geq 0.
\end{equation}

We take holomorphic linear coordinates $z_1,\dots,z_n$ with $z_j=x_j+iy_j$ on $\cn$ and suppose,
without loss of generality, the point $p$ of the lemma to be the origin and the tangent space ${\mbox{\bf T}}_0\Sigma$ of $\Sigma$ at the origin to be the real axis $\bR^n_{x_1\dots x_n}$ of $\cn$. Thus near $0$ the surface $\Sigma$ is the graph $y=\vp(x)$ for a smooth $\bR^n$--valued function $\vp$ defined on a neighborhood of $0\in\bR^n_{x_1\dots x_n}$ that satisfies $\vp(0)=0$ and $d\vp(0)=0$. There are positive constants $r$ and $k$ such that $\vp$ is defined on $D_r=\{x\in\bR^n_{x_1\dots x_n}:|x|<r\}$ and satisfies there the inequality 
$|\vp(x)|<k|x|^2$ and the derivatives of the coordinates of $\vp$ satisfy $|\partial \vp_j/\partial x_\nu|\leq k|x|$. 

Fix an $\eta$ such that $0<\eta<r$.

Define the surface $S$ to be the graph 
$$S=\{x+i\chi_\delta(|x|)\vp(x):x\in 
D_r\}$$
in which $\delta$ is chosen to satisfy $0<\delta<\min\big\{\tfrac{\eta}{2},\tfrac{r}{4}\big\}$. 
Thus the part of $S$ that lies over the annular domain $A_{2\delta;r}=\{x\in\bR^n_{x_1\dots x_n}:2\delta<|x|<r\}$ coincides with a domain in $\Sigma$. The part of $S$ that lies over the disc $\overline{D_{\delta}}$ lies in $\mbox{\bf T}_p(\Sigma)$. 
If $\Sigma_0$ is the union of $S$ and the surface obtained by excising from $\Sigma$ the graph of $\vp$ over $D_r$, 
then $\Sigma_0$ is the surface whose existence the lemma asserts.

We have to verify that if $\Sigma$ is totally real, then so is $\Sigma_0$ provided $\delta$ and $\eta$ are chosen properly, and establish the asserted approximation. For this we will show that if $\delta$ is sufficiently small, then the $\bR^n$--valued map $\chi_\delta\vp$ satisfies a Lipschitz condition with Lipschitz constant less than one on the set $\overline{D_{2\delta}}$.

On the disc $\oD_{2\delta}$ we have
$$\left|{\partial \over \partial x_\nu}\bigl(\chi_\delta(|x|)\varphi_j(x)\bigr)\right| \leq 2\mu k|x|^2/\delta \leq 8\mu k\delta.$$
On the annular domain $\overline A_{2\delta;\eta}$ the function $\chi_\delta$ is identically one so we have
$$\left|{\partial \over \partial x_\nu}\bigl(\chi_\delta(|x|)\varphi_j(x)\bigr)\right| =
\left |\partial \phi_j/\partial x_\nu \right| \leq k\eta.$$
Thus on the disc $\oD_\eta$, the gradient of $\chi_\delta \varphi_j$ satisfies
$$\| \grad \chi_\delta \varphi_j \| \leq \max\{8 \sqrt{n}\mu k\delta, \sqrt{n} k\eta\}.$$
It follows that the $\bR^n$-valued map $\chi_\delta \varphi$ satisfies a Lipschitz condition with Lipschitz constant $\max \{8n\mu  k \delta, n k \eta\}$.

Thus if $\eta< 1/ nk$ and $\delta< 1/ 8n \mu k$, then $\chi_\delta \varphi$ satisfies a Lipschitz condition with constant less than one on $\oD_\eta$.  The discussion on page 59 of \cite{Stout:2007} shows that the part $\Sigma_{\eta}$ of $\Sigma_0$ over $\oD_{\eta}$ is totally real, and Theorem~1.6.9 of the same source shows that $\sP(\Sigma_{\eta})=\sC(\Sigma_{\eta})$.

The lemma is proved.
\smallskip

For flattening functions there is the following lemma.

\begin{lemma}
\label{6.18.14.i} Let $\Sigma$ be a compact smooth $n$-dimensional submanifold of $\bC^n$, let $p\in \Sigma$, and let $\Sigma$ be totally real near $p$.  Suppose there is a smooth real-valued function $g$ on $\Sigma$ such that each level set of $g$ is polynomially convex and $g$ assumes a strict maximum at $p$.
Then for some $\varepsilon>0$, there is a smooth real-valued function $g_0$ on the manifold $\Sigma_0$ given in Lemma~\ref{4.16.14.i}
such that each level set of $g_0$ is polynomially convex and $g_0$ assumes its maximum on a neighborhood of $p$ in $\Sigma_0$. Moreover, if the level sets of $g$ satisfy $\sP(g^{-1}(t))=\sC(g^{-1}(t))$, then the level sets of $g_0$ have the same property.  
Also the range of $g_0$ coincides with the range of $g$.
\end{lemma}

Of course the same result also holds with maximum replaced throughout by minimum, and $g$ can be modified at its maximum and minimum simultaneously.

\noindent{\bf Proof.} Take $\eta>0$ as in Lemma~\ref{4.16.14.i}.
Without loss of generality $g(p)=1$.
  Since the maximum of $g$ at $p$ is strict, there is an $\alpha>0$ such that $g<1-\alpha$  on $\Sigma\setminus \bB_n(p,\eta)$.  Now choose $\varepsilon>0$ as in Lemma~\ref{4.16.14.i} with $3\varepsilon<\eta$ and such that in addition $g>1-\alpha$ on $\Sigma\cap \bB_n(p, 3\varepsilon)$.  Choose a smooth function $\psi$ on $\Sigma_0$ with \vadjust {\kern 3pt} 
$0\leq \psi\leq 1$ such that the set $\{\psi=1\}$ contains a neighborhood of $\Sigma_0\cap \overline{\bB_n(p,2\varepsilon)}$ and such that $\psi$ vanishes identically outside of $\Sigma_0\cap \bB_n(p, 3\varepsilon)$.  The function $g$ is defined on $\Sigma_0\setminus \bB_n(p, 2\varepsilon)=\Sigma\setminus \bB_n(p, 2\varepsilon)$ and if we regard $(1-\psi)g$ as zero on $\Sigma_0\cap \bB_n(p, 2\varepsilon)$, then $(1-\psi)g$ becomes a well-defined smooth function on $\Sigma_0$.  Define the function $g_0$ on $\Sigma_0$ by $g_0=\psi + (1-\psi)g$. This function is smooth and assumes its maximum value of one on a set containing $\bB_n(p, 2\varepsilon)\cap\Sigma_0$.

To see that each level set of $g_0$ is polynomially convex, first note that on $\Sigma_0\setminus \bB_n(p,\eta)$ we have $g_0=g<1-\alpha$, so each level set $\{g_0=c\}$ with $|c|\geq 1-\alpha$ is contained in $\Sigma_0\cap \overline{\bB_n(p,\eta)}$ and hence is polynomially convex by Lemma~\ref{4.16.14.i}.  
Note now that $g_0$ is identically 1 on $\Sigma_0\setminus \bB_n(p,2\varepsilon)$ and that $g_0\geq g$ on $\Sigma_0\setminus \bB_n(p, 2\varepsilon)=\Sigma\setminus \bB_n(p, 2\varepsilon)$, so if $x$ is such that $g_0(x)< 1-\alpha$, then $g(x)< 1-\alpha$ also, whence $x$ is in $\Sigma\setminus \bB_n(p, 3\varepsilon)$ and $g_0(x)=g(x)$.  Thus for each $c$ with $|c|<1-\alpha$, the level set $\{g_0=c\}$ coincides  with the set $\{g=c\}$ and so is polynomially convex by hypothesis.

The final two statements of the lemma are clear.

The lemma is proved.

\section{Attaching  Tubes.}

In this section we detail the  process of attaching tubes to surfaces.

\begin{lemma}
\label{6.22.14.cxii}
 Let $\Sigma_1$ and $\Sigma_2$ be smooth surfaces in $\bC^2$, let $a,b\in \bR$ satisfy $a<b$, and let
\begin{equation}
\label{6.23.14.i}
\bR^2=\{{(y_1,y_2)}\in \bC^2:  y_1, y_2\in \bR\}.
\end{equation}
Suppose \hfil\break
(i) $\Sigma_1\subset \{(z_1,z_2)\in \bC^2: \Re z_1\leq a\}$\hfil\break
(ii) $\Sigma_2\subset \{(z_1,z_2)\in \bC^2: \Re z_1\geq b\}$\hfil\break
(iii) $\Sigma_1$ contains a disc of radius $2r_0$ in the plane $(a,0)+i\bR^2$ centered at $(a,0)$\hfil\break
(iv) $\Sigma_2$ contains a disc of radius $2r_0$ in the plane $(b,0)+i\bR^2$ centered at $(b,0)$.

\noindent Let $\Sigma'_1$ and $\Sigma'_2$ be obtained from $\Sigma_1$ and $\Sigma_2$ respectively by excising the discs of radius $r_0$ about $(a,0)$ and $(b,0)$.  Then there is a positive function $\vr$ on the interval $[a,b]$ with $\vr(a)=\vr(b)=r_0$ such that if $S_\vr$ to is the tube given by
\begin{equation}
\label{6.24.14.i}
S_\vr=\{(t,0)+i\vr(t)(\cos\theta,\sin\theta):t\in [a,b], \theta\in [0,2\pi]\},
\end{equation}
the set $\Sigma'_1\cup\Sigma'_2\cup S_\vr$ is a smooth surface in $\bC^2$ homeomorphic to $\Sigma_1\#\Sigma_2$.
\end{lemma}

Each cross section $C_t=\{(t,0)+i\vr(t)(\cos\theta,\sin\theta): \theta\in [0,2\pi]\}$ $(a<t<b)$ of the tube is a circle lying in a totally real plane and hence is polynomially convex and satisfies $\sP(C_t)=\sC(C_t)$.

Note also that if $f$ is a continuous function on $\Sigma'_1\cup\Sigma'_2\cup S_\vr$ that is smooth on $\Sigma'_1$ and $\Sigma'_2$, is constant on each of $\Sigma'_1 \cap [(a,0)+i\bR^2]$ and $\Sigma'_2 \cap [(b,0)+i\bR^2]$, is constant on each cross section $C_t$, and increases linearly in $t$, then $f$ is smooth on 
$\Sigma'_1\cup\Sigma'_2\cup S_\vr$ and the graph of $f$ over the tube is totally real.  (The assertion about total reality follows from the remark at the end of Section~4.)
\smallskip

\noindent{\bf Proof of the Lemma.} For the function $\vr$ it suffices to take any positive smooth function on $[a,b]$ with the properties that $\vr(a)=\vr(b)=r_0$ and that the set $E$ in the plane that is the union of the graph of $\vr$, the ray $\{(a,y):y\geq r_0\}$ and the ray $\{(b,y):y\geq r_0\}$ is a smooth curve.

\section{A Function on the Torus.}

The completion of the proof of Theorem~\ref{3.31.14.iii} depends further on the construction of certain functions on the torus and on the real projective plane. For the torus, a suitable function is given in the following lemma.
\begin{lemma}
\label{7.21.1.i}
On the torus $\bT^2$ there is a smooth real-valued function $g$ with the property that each fiber $F_t=g^{-1}(t)$ is polynomially convex and satisfies $\sP(F_t)=\sC(F_t)$.  In addition, the maximum of $g$ on $\bT^2$ is $\tfrac{3}{2}$ and is taken only at the point $(1,1)$, and the minimum is $-\tfrac{3}{2}$ and is taken only at the point $(-1,-1)$.
\end{lemma}

\noindent{\bf Proof.} We take for $g$ the function  defined by
 $$g(e^{i\vt_1},e^{i\vt_2})=
\Re e^{i\vt_1}+\tfrac{1}{2}\Re e^{i\vt_2}=\cos\vt_1+\tfrac{1}{2}\cos\vt_2,$$
 The lemma is established by showing that each fiber $F_t$
is disjoint from some circle $K_a=\{(e^{i\vt_1},e^{i\vt_2})\in\bT^2:e^{i\vt_1}=a\}$ and that no fiber $F_t$ contains any complete circle $K_b$ and then invoking a lemma from \cite{I-SK-W:2020}.

If $b=\cos\beta+i\sin\beta$, then for $F_t$ to contain $K_b$ would mean that $\cos\beta+\tfrac{1}{2}\cos\vt_2=t$ for all $\vt_2$. This does not happen.

If $t\notin[-\tfrac{3}{2},\tfrac{3}{2}]$, then $F_t$ is empty. For $t\in[-\tfrac{3}{2},\tfrac{3}{2}]$ there is $\vt_1$ such that the equation $\cos\vt_1+\tfrac{1}{2}\cos\vt_2=t$, i.e, $\cos\vt_2=2(t-\cos\vt_1)$, has no solution: No matter what $t$ may be, $\vt_1$ can be chosen so that $|2(t-\cos\vt_1)|$ is greater than one. Thus $F_t$ is disjoint from the circle $K_a$ if $a=\cos\vt_1+i\sin\vt_1$.

It now follows from Lemma 3.2 of the paper \cite{I-SK-W:2020}, which was quoted above as Lemma~\ref{2.26.15.i}, 
that each of the fibers $F_t$ is polynomially convex and satisfies $\sP(F_t)=\sC(F_t)$.

The maximum of $g$, $\tfrac{3}{2}$, is strict and occurs at $(1,1)$, and the minimum of $g$, $-\tfrac{3}{2}$, is also strict and occurs at $(-1,-1)$. 
\smallskip

The surface $\bT^2$ lies 
in the slab $W$
between the real hyperplanes $\{(z_2,z_2)\in\bC^2:\Re z_1=\pm1\}$.
\smallskip

 At $(1,1)$ the tangent space of $\bT^2$   is given by $$\mbox{\bf T}_{(1,1)}\bT^2=(1,1)+i\bR^2,$$
 at $(-1,-1)$ by
$$\mbox{\bf T}_{(-1,-1)}\bT^2=(-1,-1)+i\bR^2$$
where
$$\bR^2=\{{(y_1,y_2)}\in \bC^2:  y_1, y_2\in \bR\}.$$

\section{A Function on the Real Projective Plane.}

\begin{lemma}
\label{4.4.15.i}
 The real projective plane  admits an embedding as a smooth submanifold $\bP$ of $\bC^2$ in such a way that there is a real-valued function $g\in \sC^\infty(\bP)$  such that the graph of $g$ is totally real and each fiber $F_t = g^{-1}(t)$ is polynomially convex and satisfies $\sP(F_t) = \sC (F_t)$. 
 In addition, the maximum of $g$ on $\bP$ is $\tfrac{3}{2}$ and is taken only at the point $\bigl(\tfrac{3}{4},\tfrac{\sqrt 3}{4}\bigr)$, and the minimum is $-1$ and is taken only at the point $(-1,0)$.
\end{lemma}
\noindent{\bf Proof.} We work in $\bC^2$ with coordinates $z_1 = x_1 +ix_2$ and $z_2 = x_3 +ix_4$ and in its real hyperplane	$\bR^3_{x_1 x_2 x_3}$ with equation $x_4=0$.	In 
$\bR^3_{x_1 x_2 x_3}$ we	have	the unit sphere $\bS=\{(x_1,x_2,x_3):x^2_1 +x^2_1 +x^2_3 =1\}$.

The map $\varphi = (\varphi_1,\varphi_2) : \bC^2 \rightarrow \bC^2$ given by $\varphi_1(z) = z_1^2$ and $\varphi_2(z) = z_1z_2$ satisfies $d\varphi_1\wedge d\varphi_2= 2z_1^2\, dz_1\wedge dz_2$
so that it is locally biholomorphic off the $z_2$--axis. This implies that the restriction of $\varphi$ to $S$ is local diffeomorphic onto its range except possibly at the two points $(0, \pm 1)$. At the point $(0,1)$ the tangent plane to $\bS$ is the complex line $\l = \{(z_1,1) : z_1\in\bC\}$. As $\varphi(z_1, 1) = (z_1^2, z_1)$, which is locally diffeomorphic on $\l$ near $(0,1)$, it follows that the restriction $\varphi |\bS$ is locally diffeomorphic near the point $(0,1)$.  Similarly, $\varphi |\bS$ is locally diffeomorphic near the point $(0, -1)$. The map $\varphi$ is invariant under the involution $\varepsilon:\bC^2\rightarrow \bC^2$ given by $\varepsilon(z) = -z$, so $\varphi$ carries $\bS$ onto its image $\bP$ in $\bC^2$ as a two-to-one covering map.  The smooth submanifold $\bP$ of 
$\bC^2$ is a diffeomorphic copy of the real projective plane.

The only complex lines contained in $\bC\times \bR=\bR^3_{x_1 x_2 x_3}$ are those of the form $\bC\times \{c\}$ for a fixed $c\in \bR$, so the only complex tangents to $\bS$ occur at the points $(0,\pm 1)$.  Because the derivative of $\varphi$ is an invertible 
complex-linear transformation at all other points of $\bS$, it follows that $\bP$ has no complex tangents at points other than $\varphi(0,\pm 1)=(0,0)$.  The tangent space to $\bP$ at the point $(0,0)$ is the complex line $\{0\}\times \bC$.

Let $g$ be the smooth real-valued function on $\bC^2$ defined by 
\begin{equation}
\label{4.30.15.i}
g(z_1,z_2)=\Re(z_1+\sqrt 3 z_2)=x_1 + \sqrt 3 x_3.
\end{equation}
By the discussion in the section of preliminaries on totally real embeddings, to show total reality of the graph of $g$ over $\bP$ it suffices to check totally reality at the point over $(0,0)$.  
This one can do by applying Lemma~\ref{6.20.14.i}.
 Explicitly, notice that with holomorphic coordinates $(w_1,w_2,w_3)$ on $\bC^3$ and with $\Phi:\bC^2\rightarrow\bC^3$ the graph map given by $\Phi(z_1,z_2)=(z_1,z_2,g(z_1,z_2))$, we have that 
\begin{equation*}
\begin{split}
\Phi^*(dw_2\wedge dw_3)&=dz_2\wedge dg\\
&=-dx_1\wedge dz_2-\sqrt 3 idx_3
\wedge dx_4,
\end{split}
\end{equation*}
which does not vanish on the tangent plane $\mbox{\bf T}_{(0,0)}\bP$.

For the proof that each fiber $F_t=g^{-1}(t)$ of $g$ on $\bP$ is polynomially convex, first note that given $t\in \bR$, the function $h_t(z)=i(z_1+\sqrt 3 z_2 - t)$ is real-valued on $F_t$, and consequently to show $\sP(F_t)=\sC(F_t)$ it suffices, by Lemma~3.1, to show that each level set $L$ of $h_t$ on $F_t$ satisfies $\sP(L)=\sC(L)$.  Since each level set of $h_t$ on $F_t$ is simply a level set of the function $f(z_1,z_2)=(z_1+\sqrt 3 z_2)$, we can conclude that each fiber $F_t$ satisfies $\sP(F_t)=\sC(F_t)$ if we can show that each level set of the function $f$ is a finite set.

To show each level set of $f$ on $\bP$ is finite we show each level set of the function $f\circ \varphi$ on $\bS$ is finite.

On $\bS$ the function $f\circ\vp$ is given by
\begin{equation*}
(f\circ\vp)(z_1,z_2)=(x_1^2-x_2^2+\sqrt 3 x_1x_3)+i(2x_1x_2+\sqrt 3 x_2x_3),
\end{equation*}
so the level set 
$$
\{(z_1,z_2)\in\bS:(f\circ\vp)(z_1,z_2)=t+iw\}
$$
is given by the system of three simultaneous equations
\begin{equation}
\label{3.27.15.i}
\begin{cases}
&x_1^2-x_2^2+\sqrt 3 x_1x_3=t\\
&2x_1x_2+\sqrt 3x_2x_3 = w\\
&x_1^2+x_2^2+x_3^2=1.\\
\end{cases}
\end{equation}

Using the last equation to eliminate $x_2$ from the first two equations yields

\begin{equation}
\label{3.27.15.iii}
\begin{cases}
&x_1^2-(1-x_1^2-x_3^2)+\sqrt 3 x_1x_3=t\vadjust{\kern 4pt}\\
&\pm(2x_1+\sqrt 3 x_2)\sqrt{1-x_1^2-x_3^2}=w.\\
\end{cases}
\end{equation}
Rearranging the first equation and squaring the second lets us rewrite this system as
\begin{equation}
\label{3.24.15.i}
\begin{cases}
&2(x_1+\tfrac{\sqrt 3}{4}x_3)^2+\tfrac{5}{8}x_3^2=t+1\vadjust{\kern 4pt}\\
&(2x_1+\sqrt 3 x_3)^2(1-x_1^2-x_3^2)=w^2.\\
\end{cases}
\end{equation}
Let $u_1=x_1+\tfrac{\sqrt 3}{4}x_3$ and $u_3=x_3$. In $(u_1,u_3)$--coordinates the system (\ref{3.24.15.i}) is
\begin{equation}
\label{3.24.15.ii}
\begin{cases}
\noindent \quad\quad 2u^2_1+\tfrac{5}{8} u_3^2=t+1\vadjust{\kern 4pt}\\
\quad\quad (2u_1+\tfrac{\sqrt 3}{2}u_3)^2(1-u_1^2+\tfrac{\sqrt 3}{2}u_1u_3-\tfrac{19}{16}u_3^2)=w^2\\
\end{cases}
\end{equation}

Note that the solution set of the first equation of the system (\ref{3.24.15.ii}) is empty for $t<-1$ and that for $t=-1$ it is the single point $(0,0)$. Thus, we have only to consider the case that $t>-1$, in which case the solution set of the first equation of (\ref{3.24.15.ii}) is an ellipse $E$.

The set $L$ in $(u_1,u_3)$--space determined by the system (\ref{3.24.15.ii})
is a real-analytic variety that is a subset of $E$. As such, it is either all of $E$ or else is a discrete, and hence finite, subset of $E$. It thus suffices to show that $L$ is not all of $E$.

Let $p=p(u_1,u_3)$ and $q=q(u_1,u_3)$ be the functions on $\bR^2$ defined by the expressions on the left sides of the system (\ref{3.24.15.ii}).
 To show that $L$ is not all of $E$, it is enough to show that there is a point on $E$ where $\grad p$ and $\grad q$ are linearly independent.
 When $u_3=0$, we have $\tfrac{\partial p}{\partial u_3}=0$, so it is enough to show that $\tfrac{\partial q}{\partial u_3}\neq 0$ when $u_3=0$. Computing yields
\begin{equation*}
\begin{split}
\tfrac{\partial q}{\partial u_3}\big|_{u_3=0}&=2(2u_1)\tfrac{\sqrt 3}{2}(1-u_1^2)+(2u_1)^2\bigl(\tfrac{\sqrt 3}{2}u_1\bigr)\\
&=2\sqrt 3u_1(1-u_1^2)+2\sqrt 3u_1^3\\
&=2\sqrt 3u_1\\
&\neq 0,
\end{split}
\end{equation*} 
the latter inequality because $u_1\neq 0$ when $ u_3=0$, granted that $t>-1$ as we are assuming.

This concludes the proof that $\sP(F_t)=\sC(F_t)$ for each fiber $F_t=g^{-1}(t)$.
\smallskip

It remains to show that the maximum and the minimum of $g$ are each attained at a unique point of $\bP$. For this, it is convenient to work with the function $g\circ\vp$ on $\bS$. In real coordinates, $g\circ\vp$ is given by $$g\circ\vp(x_1,x_2,x_3)=x_1^2-x_2^2+\sqrt 3x_1x_3.$$ Eliminating $x_2$ as before reduces our problem to that of finding the extrema of the function 
$\tilde g$ given by 
$$\tilde g(x_1,x_3)=2x_1^2+x^2_3+\sqrt 3 x_1x_3-1$$
on the unit disc $$\{(x_1,x_3):x_1^2+x_3^2\leq 1\}$$
in the $(x_1,x_3)$--plane. The gradient of $\tilde g$ is 
$$\grad\tilde g(x_1,x_3)=(4x_1+\sqrt 3 x_3,2x_3+\sqrt 3 x_1),$$
which vanishes only at the point $(0,0)$. This corresponds, in complex notation, to the points $(\pm i,0)$ on $\bS$ and to $(-1,0)$ on $\bP$, where $g$ takes the value $-1$.

The extrema of $\tilde g$ on the boundary, i.e., on the circle with equation $x_1^2+x_3^2=1$, can be found by the method of Lagrange multipliers. This yields the system 
\begin{equation*}
\begin{cases}
4x_1+\sqrt 3 x_3&=\lambda x_1\\
\sqrt 3 x_1+2x_3&=\lambda x_3\\
x_1^2+x_3^2&=1
\end{cases}
\end{equation*}
The first two of these equations imply that $\lambda=1$ or $5$. The solution with $\lambda= 5$ yields the two points $(x_1,x_3)=\bigl(\pm\tfrac{\sqrt 3}{2},\pm\tfrac{1}{2}\bigr)$. We have $\tilde g\bigl(\pm\tfrac{\sqrt 3}{2},\pm\tfrac{1}{2}\bigr)=\tfrac{3}{2}$. The solution with $\lambda=1$ yields the two points $(x_1,x_3)=(\pm\tfrac{1}{2},\mp\tfrac{\sqrt 3}{2}\bigr)$  and $\tilde g\bigl(\pm\tfrac{1}{2},\mp\tfrac{\sqrt 3}{2}\bigr)=-\tfrac{1}{2}$.

Putting this all together, we conclude that the maximum of $g$ occurs at the point $\vp\bigl(\pm\tfrac{\sqrt 3}{2},\pm\tfrac{1}{2}\bigr)$ and is $\tfrac{3}{2}$ and that the minimum of $g$ occurs at $(-1,0)$ and is $-1$.

\section{Conclusion of the proof of Theorem~\ref{3.31.14.iii}.}

In this section we conclude the proof of Theorem \ref{3.31.14.iii} 
in the case of closed surfaces by showing that every compact closed surface is diffeomorphic to a surface that satisfies the hypotheses of Lemma~\ref{2.27.15.i}.  
In $\bC^2$ we work with coordinates
$z_1 = x_1 +ix_2$ and $z_2 = x_3 +ix_4$.

We first treat the case of the sphere.  In $\bC^2$  we have the sphere 
$\bS=\{(x_1,x_2,x_3):x^2_1 +x^2_1 +x^2_3 =1, x_4=0\}$ and the function $x_1+2$ on $\bS$.  Applying the flattening process of Lemmas~\ref{4.16.14.i} and~\ref{6.18.14.i} to the sphere $\bS$ and function  $x_1+2$ at the point $p=(-1,0)$ yields a surface and function satisfying the hypotheses of Lemma~\ref{2.27.15.i}. The only complex tangents to the flattened sphere occur at the points $(0,1)$ and $(0,-1)$, and at those points $d(x_1+2)\neq 0$. Thus by Theorem~\ref{2.27.15.i}, the theorem holds in the case of the sphere.

Next we treat the case of the connected sum of a finite number, say $m$, of tori.  The standard torus $\bT^2$ contains the points $(-1,-1)$ and $(1,1)$ and supports the function $g$ constructed in Section 9 above.

We denote by $T_0$ the torus obtained from $\bT^2$ by applying the flattening process of Lemma~\ref{4.16.14.i} to the torus $\bT^2$ at both the points $(1,1)$ and $(-1,-1)$. We then apply the flattening process of Lemma~\ref{6.18.14.i} to obtain from the function $g$ constructed in Section 9 a new function $g_0$ on the torus $T_0$ 
also flattened at both the points $(1,1)$ and $(-1,-1)$. 
Note that the range of $g_0$ is the interval $\bigl[-\tfrac{3}{2},\tfrac{3}{2}\bigr]$.

Introduce the translation $\tau:\bC^2\rightarrow\bC^2$ given by
\begin{equation}
\label{6.23.14.iii}
\tau(z_1,z_2)=(z_1+5,z_2+2).
\end{equation}
Thus,
$\tau(-1,-1)=(4,1),\tau(1,1)=(6,3), \tau^2(-1,-1)=(9,3), \tau^2(1,1)=(11,5),\ldots$. The point to be observed 
is that the second coordinate of $\tau^{n+1}(-1,-1)$ is the same as second coordinate of $\tau^n(1,1)$. This is evident since the second coordinates of $(-1,-1)$ and $(1,1)$ differ by two. 
Also note that the first coordinate of $\tau^{n+1}(-1,-1)$ is strictly greater than the first coordinate of $\tau^n(1,1)$.

For $n=0,\ldots, m-1$,
let $\tau^n(1,1)=(\xi_n,\eta_n)$, $\tau^n(-1,-1)=(\xi_n^-,\eta_n^-)$, and let $T_n$ be the torus $\tau^n(T_0)$. (By definition $\tau^0$ is the identity map.)
The torus $T_n$ is contained in the strip
$W_n=\{(z_1,z_2):\xi_n^-\leq \Re z_1\leq \xi_n\}$.
Note that the strips $W_0, W_1, \ldots, W_{m-1}$ are pairwise disjoint.  The straight line interval $\ell_n$ connecting $(\xi_n,\eta_n)$ to $(\xi^-_{n+1},\eta_{n+1}^-)$ lies in the line 
$$L_n=\bigl\{\bigl((1-s)\xi_n+s\xi^-_{n+1},\eta_n\bigr):s\in\bR\bigr\}=\big\{\bigl((1-s)\xi_n+s\xi^-_{n+1},\eta^-_{n+1}\big):s\in\bR\big\}.$$

Let $T_n'$ be the torus $T_n$ from which suitable small discs centered at $\tau^n(-1,-1)$ and at $\tau^n(1,1)$ have been excised, the discs chosen so that 
the lemma on tubes, Lemma~\ref{6.22.14.cxii}, yields tubes $S_{\vr_j}$ of the form (\ref{6.24.14.i}), 
$j=0,\ldots,m-2$, $S_{\vr_j}$ centered along the interval 
$\ell_j$, such that the surface $\mathcal T$ given by
$$\mathcal T=T_0'\cup\dots\cup T_{m-1}'\cup S_{\vr_1}\cup\cdots \cup S_{\vr_{m-2}}$$
is a smooth surface, which is, topologically, the connected sum of $m$ tori.

We define on the surface $\mathcal T$ a smooth function $f$ by the conditions that on $T_k$, $f$ is to coincide with the function $g_0\circ\tau^{-k}+4k$. With this prescription, the range of $f$ on $T_{k+1}$ is to the right of the range of $f$ on $T_k$ and is at distance one from it. Accordingly, we can require that the function $f$ on $S_{\vr_k}$ be constant on the cross sections of $S_{\vr_k}$ and increase linearly from $T_k$ to $T_{k+1}$.

The surface $\mathcal T$ and the function
$f+{5\over 2}$ satisfy the hypotheses of Lemma~\ref{2.27.15.i} above.  
Because $\bT^2$ is totally real and the flattening process preserves total reality, the only complex tangents to $\mathcal T$ occur at points of the tubes, and along the tubes $d(f+{5\over 2})\neq 0$ because $f$ increases linearly along the tubes. Consequently, appealing to Lemma~\ref{2.27.15.i}  completes the proof of the theorem in the case of orientable surfaces.
\smallskip

To conclude the proof of the theorem, we treat the case of nonorientable surfaces. 

Consider a surface $\mathcal S$ that is the connected sum of $m$ copies of the real projective plane.

We have the diffeomorphic copy $\bP$ of the real projective plane and the function $g$ on it given by Lemma \ref{4.4.15.i}. The function $g$ is real-valued, has range $[-1,\tfrac{3}{2}]$,  and has polynomially convex level sets. By construction $g$ takes the value $-1$ at the point $(-1,0)$ and nowhere else, and it takes the value $\tfrac{3}{2}$ at the point $\bigl(\tfrac{3}{4},\tfrac{\sqrt 3}{4}\bigr)$ and nowhere else. The graph of $g$ is totally real.

Introduce the surface $\bP_0$ obtained from $\bP$ by applying the flattening process of Section 7 at the points $(-1,0)$ and $\bigl(\tfrac{3}{4},\tfrac{\sqrt 3}{4}\bigr)$.
Concerning the surface $\bP_0$ we note that it coincides with $\bP$ outside two small totally real discs, one containing $(-1,0)$ the other $\bigl(\tfrac{3}{4},\tfrac{\sqrt 3}{4}\bigr)$. Consequently, it is totally real except at the point $(0,0)$. 

Let $g_0$ be the function obtained from the function $g$ on $\bP$ by flattening $g$ at the points points $(-1,0)$ and $\bigl(\tfrac{3}{4},\tfrac{\sqrt 3}{4}\bigr)$ in accordance with the flattening process of the same section. The function $g_0$ has range $[-1,\tfrac{3}{2}]$, the graph of $g_0$ is totally real, and the level sets of $g_0$ are polynomially convex.

The function $g$ ({\it{not}} $g_0$) can, in fact, be defined on all of $\bC^2$ by the formula 
\begin{equation*}
g(z_1,z_2)= x_1+\sqrt 3 x_3.
\end{equation*}
(See equation \ref{4.30.15.i} in the proof of Lemma~\ref{4.4.15.i}.)
Then the surface $\bP_0$ is contained in the strip $W$ in $\bC^2$ given by
$$W=\{(z_1,z_2)\in\bC^2:-1\leq g(z_1,z_2)\leq \tfrac{3}{2}\}.$$
The points $(-1,0)$ and $\bigl(\tfrac{3}{4},\tfrac{\sqrt 3}{4}\bigr)$ lie in the real hyperplanes that bound $W$, one in one of the boundary components, one in the other.

In addition to $\bP_0$, introduce its inversion $\tilde \bP_0=\varrho(\bP_0)$ with $\varrho:\bC^2\rightarrow\bC^2$ the inversion defined by
$$\varrho(z_1,z_2)=(-z_1,-z_2).$$
That is $\tilde \bP_0$ is the surface obtained by reflecting $\bP_0$ through the origin.  Note that $\varrho^2=\varrho\circ\varrho$ is the identity map and
that the surface $\bP_0$ is contained in the strip $\varrho(W)$ given by
$$\varrho(W)=\{(z_1,z_2)\in\bC^2:-\tfrac{3}{2}\leq g(z_1,z_2)\leq 1.$$

On the surface $\tilde \bP_0$ define the function $\tilde g_0$ by the condition that
$$\tilde g_0(z_1,z_2)=\, -g_0\circ \varrho.$$
The function  $\tilde g_0$ takes its maximum of $1$ on a neighborhood of the point $(1,0)$ and its minimum of $-\tfrac{3}{2}$ on a neighborhood of the point $\bigl(-\tfrac{3}{4},-\tfrac{\sqrt 3}{4}\bigr)$, the graph of $\tilde g_0$ is totally real, and the level sets of $\tilde g_0$ are polynomially convex.

Note that the tangent space to $\tilde \bP_0$ at the point $(1,0)$ (regarded as a real vector subspace of $\bC^2$) coincides with the tangent space to $\bP_0$ at the point $(-1,0)$, and the same holds for the tangent spaces to $\tilde \bP_0$ at the point $\bigl(-\tfrac{3}{4},-\tfrac{\sqrt 3}{4}\bigr)$ and to $\bP_0$ at the point $\bigl(\tfrac{3}{4},\tfrac{\sqrt 3}{4}\bigr)$.  Because all of these tangent spaces are 2-dimensional totally real subspaces of $\bC^2$ and the null space of the real-linear function $g$ is of real dimension 3 and contains each of these tangent spaces, we can choose vectors $\vec u$ and $\vec v$ such that $i\vec u$ lies in the tangent space to $\bP_0$ at $(-1,0)$ while $g(\vec u)\neq 0$ and $i\vec v$ lies in the tangent space to $P_0$ at $\bigl(\tfrac{3}{4},\tfrac{\sqrt 3}{4}\bigr)$ while $g(\vec v)\neq 0$.  By rescaling $\vec u$ and $\vec v$, we may arrange to have $g(\vec u)=g(\vec v)=10$.

{Let $\tau_1,\ldots, \tau_{m-1}:\bC^2\rightarrow \bC^2$ be translation operators
$$\tau_k(z_1,z_2)=(z_1,z_2)+\vec a_k$$
and define projective planes $\bP_1,\ldots \bP_{m-1}$ by setting $\bP_k=\tau_k(\bP_0)$ for $k$ even, and  $\bP_k=\tau_k(\tilde \bP_0)$ for $k$ odd.  Choose the $\vec a_k$ defining the translations such that for $k$ even
$$\tau_{k+1}\bigl(-\tfrac{3}{4},-\tfrac{\sqrt 3}{4}\bigr)-\tau_k\bigl(\tfrac{3}{4},\tfrac{\sqrt 3}{4}\bigr)=\vec v$$
and for $k$ odd
$$\tau_{k+1}(-1,0)-\tau_k(1,0)=\vec u.$$}

{For each $k$, the range of the real-linear function $g$ on $\bP_{k+1}$ lies strictly to the right of the range of $g$ on $\bP_k$.  Thus the $\bP_k$ lie in pairwise disjoint strips defined by inequalities on $g$.  Deleting small discs centered at the points on each $\bP_k$ corresponding to the points $(-1,0)$ and 
$\bigl(\tfrac{3}{4},\tfrac{\sqrt 3}{4}\bigr)$ on $\bP_0$ to form surfaces $\bP'_0,\ldots, 
\bP'_{m-1}$ and attaching thin tubes centered along the line segments from 
$\tau_k\bigl(\tfrac{3}{4},\tfrac{\sqrt 3}{4}\bigr)$ to $\tau_{k+1}\bigl(-\tfrac{3}{4},-\tfrac{\sqrt 3}{4}\bigr)$ when $k$ is even and from $\tau_k(1,0)$ to $\tau_{k+1}(-1,0)$ when $k$ is odd yields a surface that is topologically the connected sum of $m$ copies of the projective plane.}

 We now apply the lemma of Section~10 above on attaching tubes.  Consider the case of $k$ odd.  Because the tangent space to $\bP_k$ at $\tau_k(-1,0)$ and the tangent space to $\bP_{k+1}$ at $\tau_{k+1}(-1,0)$ coincide and are totally real, and $i[\tau_{k+1}(-1,0)-\tau_k(1,0)]=i\vec u$ lies in that totally real space, which we will denote by $\Pi$, there is a complex-linear automorphism $\Phi:\bC^2\rightarrow\bC^2$ that takes $\Pi$  to $i\bR^2$ and takes $i[\tau_{k+1}(-1,0)-\tau_k(1,0)]$ to the vector $(1,0)$.  Allowing $\Phi$ to be 
 complex-affine rather than complex-linear  
allows further for $\Phi$ to send $\tau_{k}(1,0)$ and $\tau_{k+1}(-1,0)$ to the points $(0,0)$ and $(1,0)$, respectively.  Then setting $\Sigma_1=\Phi(P_k)$ and $\Sigma_2=\Phi(P_{k+1})$, we are exactly in the situation where Lemma~\ref{6.22.14.cxii} on attaching tubes applies to yield a tube $S_\varrho$ connecting $\Sigma'_1$ and $\Sigma'_2$.  The tube $S_{\varrho_k}=\Phi
^{-1}(S_\varrho)$ then connects $\bP'_k$ and $\bP'_{k+1}$.  The case of $k$ even is handled in the same manner.  The surface 
$${\mathcal P}= \bP'_0\cup\dots\cup\bP'_{m-1}\cup S_{\varrho_0}\cup\dots\cup S_{\varrho_{m-2}}$$
is a smooth surface that is topologically the connected sum of $m$ copies of the projective plane.

We define on the surface $\mathcal P$ a smooth function $f$ by the condition that for $k$ even, $f$ is to coincide on $\bP_k$ with the function $g_0\circ \tau_k^{-1} + 10k$, and for $k$ odd  $f$ is to coincide on $\bP_k$ with $\tilde g_0\circ \tau_k^{-1} +10k$.
With this prescription, the range of $f$ on $\bP_{k+1}$ is strictly to the right of the range of $f$ on $\bP_k$. Accordingly, we can require that the function $f$ on $S_{\vr_k}$ be constant on the cross sections of $S_{\vr_k}$ and increase linearly from $\bP_k$ to $\bP_{k+1}$.

The surface $\mathcal P$ and the function
 {$f+2$} satisfy the hypotheses of Lemma~\ref{2.27.15.i} above. 
We have already essentially observed that the part of the graph of $f+2$ over $\bP'_0\cup\dots\cup\bP'_{m-1}$ is totally real, and the part of the graph over the tubes is totally real because $d(f+2)\neq0$ there.  Consequently, appealing to Lemma~\ref{2.27.15.i}  completes the proof of the theorem in the case of nonorientable surfaces.

This concludes the proof of the theorem for closed surfaces.  
The case of surfaces with boundary is treated in the next section.

\medskip

It should be observed that for most closed surfaces, the main theorem, Theorem~\ref{3.31.14.iii}, can be established without recourse to Theorem~\ref{9.3.14.i} with its rather involved proof. For this recall that every closed surface other than the sphere, the real projective plane, or the Klein bottle is of the form $\bT^2\#\Sigma$ for another closed surface $\Sigma$. The torus contains a compact set $E$ with $\what E\setminus E$ nonempty but devoid of analytic discs. Using the constructions given above, one finds with essentially no further work a smooth function $f$ on $\bT^2\#\Sigma$ whose graph, $\Gamma_f$, a subset of $\bC^3$, is a closed surface diffeomorphic to $\bT^2\#\Sigma$ of the sort whose existence is asserted by Theorem~\ref{3.31.14.iii}.

If in constructing embeddings of surfaces in complex Euclidean spaces so as to have hull without analytic structure, one is willing to settle for embeddings in $\C^8$ rather than $\C^3$, then a simpler proof is possible in all cases.
In the paper \cite{I-SK-W:2020} the authors use the fact that every smooth manifold of dimension $d$, $d\geq 3$,  contains a two-dimensional torus together with the result of Alexander \cite{Alexander:1998} to show that every compact manifold $S$ of dimension $d$  at least three can be embedded as a smooth submanifold $\Sigma$ of $\bC^{2d+4}$ in such a way that $\what\Sigma\neq \Sigma$ but $\what\Sigma$ contains no analytic disc.  Because every neighborhood of the diagonal in $\bT^2$ contains a set $E$ such that $\widehat E \neq E$ but $\widehat E$ contains no analytic disc (as noted at the end of the introduction), one can replace the torus in that argument by an annulus.  Since every surface contains an annulus, one then sees that the argument in \cite{I-SK-W:2020} can be repeated, essentially without change, to embed each compact surface $S$ in $\bC^8$ as a submanifold $\Sigma$ with the property that $\what\Sigma$ contains no analytic discs although $\what\Sigma$ is strictly larger than $\Sigma$.

\section{Surfaces with Boundary.}

We now show that our main result obtains for compact surfaces with boundary as a consequence of the result for closed surfaces.

Let $S$ be a smooth compact surface with nonempty boundary. Thus $S$ is obtained from a smooth compact surface $\tilde S$ without boundary by excising a finite number, say $\nu$, of open discs with smooth boundary.  Let $\Sigma$ be a smooth submanifold of $\cn$ that is diffeomorphic to  $\tilde S$. Assume 
that $\Sigma$ contains a closed subset $E$ such that $\what E \neq E$, that $\what E$ contains no analytic disc, that $\what E\cap\Sigma=E$, and that $\what \Sigma=\Sigma\cup \what E$. Finally, assume that $\Sigma\setminus E$ contains a totally real disc $U$. 
The existence of a surface $\Sigma$ satisfying all these assumptions has been proven above.
With this setup, we have the following fact:

\begin{lemma}
\label{5.4.14.iii}
The surface $\Sigma$ contains a smooth copy $\Sigma^*$ of $S$ that contains $E$ and has the property that $\widehat{\Sigma^*}=\Sigma^*\cup \what E$. 
\end{lemma}

Thus, $\what{\Sigma^*}\setminus \Sigma^*$ is nonempty but contains no analytic disc.
\smallskip

\noindent{\bf Proof.}  Let $\Sigma^*$ be the surface obtained from
the surface $\Sigma$ by excising $\nu$ open discs $D_1,\dots,D_\nu$ whose closures are contained in $U$, that have smooth boundaries, that have mutually disjoint closures, and that are small enough that $\sP(\bar D_j)=\sC(\bar D_j)$ for each $j$. This surface is diffeomorphic to the initially given surface $S$. We will show that  $\widehat{\Sigma^*}=\Sigma^*\cup \what E$. 

As $\Sigma^*\subset\Sigma$, we have $\what\Sigma^*\subset \what \Sigma$. We have to show that no point of any of the discs $D_j$ is in $\what{\Sigma^*}$ and that $\what\Sigma_0''$  contains $\what E$.

That $\what{\Sigma^*}$ contains $\what E$ is immediate, for $\Sigma^*$ contains $E$.

That no point of a $D_j$ lies in $\what{\Sigma^*}$ seems to lie deeper. The fact is that each point of each $D_j$ is a peak point for the algebra $\sP(
\what\Sigma)$. This follows from a general result in the paper \cite[Lemma~5.6]{Izzo:2003} by one of the present authors. It is also a consequence of an earlier, less general,  
result of Freeman \cite[Theorem~2.3]{Freeman:1966}. Alternatively, it is an immediate consequence of the principal result of the paper \cite{Allan:1971} of Allan.

The lemma is proved.

\section{Polynomially Convex Embeddings.}

 Each compact {\it{orientable}} surface embeds in $\bR^3\subset \bC^3$.  As a compact subset of $\bR^3$, the embedded surface is a polynomially convex subset of $\bC^3$. For homological reasons,$^($\footnote{For orientable manifolds, this is a result of Browder \cite{Browder:1961}; an exposition of the full result is in \cite[Section 2.3]{Stout:2007}.}$^)$  closed surfaces cannot be realized as polynomially convex subsets of $\bC^2$. The methods we have used suffice to show that nonorientable surfaces also can be embedded in $\bC^3$ as polynomially convex sets, though their nonorientability precludes their being embedded in $\bR^3$. To get such a polynomially convex embedding of a nonorientable surface
  $S$, it suffices to take the graph of a smooth real-valued function $h$ on a smoothly embedded copy of $S$ in $\bC^2$ such that each level set of $h$ is polynomially convex. Such functions have been constructed above. Thus we have that {\it{every closed surface embeds in $\bC^3$ as a smooth polynomially convex surface.}}
 Note also that our arguments show that the embedded surface in $\bC^3$ can always be taken to be totally real.

Franc Forstneri\v c has drawn our attention to the papers
\cite{Forstneric:1994}, \cite{Forstneric-Rosay:1993}, and  \cite{Low-Wold:2009}, which contain much deeper results on polynomially convex totally real embeddings of manifolds. In particular, in \cite{Forstneric-Rosay:1993}, Forstneri\v c and Rosay show that {\it{for a smooth compact surface $\Sigma$, 
 the set of embeddings $\vp$ of $\Sigma$ in $\bC^m, m
\geq 3$, such that $\vp(M)$ is polynomially convex and totally real is dense in the space of all embeddings of $\Sigma$ in $\bC^m.$}} In fact, the result of \cite{
Forstneric-Rosay:1993} is much stronger than this; we refer to the original source for the precise formulation.

\section{Remarks and Open Questions.}

1. Our examples are {\it smooth} manifolds, but they are not real-analytic. Conjecture: {\it{If $\bM$ is a compact real-analytic submanifold of $\cn$ such that the set $\what\bM\setminus\bM$ is nonempty, then $\what{\bM}$ contains an analytic disc.}} Note that, according to a result of Diederich and Forn\ae ss,
 \cite{Diederich-Fornaess:1978} or \cite[pp. 334-335]{Stout:2007}, $\bM$ itself can contain no analytic disc.
\smallskip

2. Are there examples in $\bC^2$? That is to say, is there a compact surface $S$ in $\bC^2$ such that $\what S$ is bigger than $S$ but contains no analytic disc?
Note that by a result of Alexander \cite{Alexander:1996}, in contrast to the examples above, no such surface $S$ could be totally real.
\smallskip

3. Does every smooth closed surface in $\bC^2
$ contain a compact set $E$ such that $\what E\setminus E$ is not empty but contains no analytic disc? There are very simple examples of smooth compact surfaces with boundary in $\bC^2$ that do not contain such sets, e.g., the polynomially convex annulus $\{(z_1,z_2)\in\bC^2:z_1z_2=1\ \mbox{and}\ 1\leq |z_1|\leq 2\}$.

4. For a fixed positive integer $n$, what is the smallest integer $p$ with the property that every compact $n$--dimensional smooth manifold can be smoothly embedded in $\bC^p$ as a polynomially convex submanifold?

5. For a fixed positive integer $n$, what is the smallest integer $p$ with the property that every compact $n$--dimensional smooth manifold can be embedded as a smooth manifold $\Sigma$ in $\bC^p$ so that $\widehat \Sigma\setminus \Sigma$ is nonempty but contains no analytic disc.

6. Purvi Gupta (private communication) has raised the question of whether each compact surface can be smoothly embedded in $\C^3$ so as to be rationally convex in addition to having the properties in Theorem~\ref{3.31.14.iii}.  This question is open in general, and we do not know whether the surfaces constructed above using Theorem~\ref{9.3.14.i} are rationally convex.  However, for every closed surface other than the sphere, the real projective plane, or the Klein bottle, the construction discussed in the penultimate paragraph of Section~11 yields a rationally convex embedding as a surface $X$ with $R(X)=C(X)$.  (Here $R(X)$ denotes the uniform closure of the rational functions  holomorphic on (a neighborhood of) $X$.)  Reason:  The surface $X$ is the graph in $\C^3$ of a smooth real-valued function $f$ on a surface $\bT^2\#\Sigma$ obtained by attaching a surface $\Sigma$ via a tube to the standard torus $\bT^2$.  The function $f$ can be taken so that the zero set of $f$ is the set $E\subset\bT^2$ of Alexander with polynomially convex hull without analytic discs and each of the other level sets $L$ of $f$ is polynomially convex and satisfies $P(L)=C(L)$.  To show that $R(X)=C(X)$, it suffices, by the Bishop antisymmetric decomposition, to show that 
$R(\bT^2\#\Sigma)|L$ is dense in $C(L)$ for each level set $L$ of the real-valued function $f$.  The required density is immediate for each level set other that the zero set $E$.  The set $E$ is contained in the standard torus $\bT^2$.  The way $X$ was constructed, the coordinate function $z_1$ never vanishes on $X$.  Thus $1/z_1$ is in $R(X)$, and on $E$ we have $1/z_1=\overline z_1$.  Thus our problem further reduces to proving approximation on the intersection of $E$ with each level set of $z_1$.  But each such intersection is a proper subset of a circle in the $z_2$ variable and so 
the set of polynomials, and hence $R(X)$, is dense there.  We conclude that 
$R(X)=C(X)$, as desired.

\bibliographystyle{plain}
\markboth{Bibliography}{Bibliography}
\bibliography{Izzo.Stout.3}

\medskip

\end{document}